\documentclass[smallextended]{svjour3}       
\usepackage[bottom=4cm, right=4cm, left=4cm, top=4cm]{geometry}

\usepackage{latexsym}
\usepackage{amssymb}
\usepackage{amsmath}
\usepackage[mathscr]{eucal}
\usepackage{graphicx}
\usepackage{hyperref}
\usepackage{caption}
\usepackage{subcaption}

\renewcommand{\qed}{\hfill{\ \ \rule{2mm}{2mm}} \vspace{0.2in}}

\newcommand{\ind}{1\hspace{-2.3mm}{1}}

\begin{document}



\title{Extremal Planar Matchings of Inhomogenous Random Bipartite Graphs}
\author{ \textbf{Ghurumuruhan Ganesan}
\thanks{E-Mail: \texttt{gganesan82@gmail.com} } \\
\ \\
University of Bristol, UK}
\date{}
\maketitle

\begin{abstract}
In this paper we study maximum size and minimum weight planar matchings of inhomogenous random bipartite graphs. Our motivation for this study comes from efficient usage of cross edges in relay networks for overall improvement in network performance. We first consider Bernoulli planar matchings with a constraint on the edge length and obtain deviation estimates for the maximum size of a planar matching. We then equip each edge of the complete bipartite graph with a positive random weight and obtain bounds on the minimum weight of a planar matching containing a given number of edges. We also use segmentation and martingale methods to obtain~\(L^2-\)convergence of the minimum weight, appropriately scaled and centred.

\vspace{0.1in} \noindent \textbf{Key words:} Length Constrained Planar Matchings, Minimum Weight Planar Matchings, Weight Saturation.

\vspace{0.1in} \noindent \textbf{AMS 2000 Subject Classification:} Primary: 60C05;
\end{abstract}

\bigskip

\renewcommand{\theequation}{\thesection.\arabic{equation}}
\setcounter{equation}{0}
\section{Introduction} \label{intro}
A graph matching, i.e., a set of vertex disjoint edges are important objects of  study and particular, matchings in bipartite graphs have well-known applications in  various fields including engineering and economics. In this paper, we focus a specific type of matchings called planar matchings which include an additional condition that the edges of the matching do not intersect each other. Planar matchings in  graphs are used to determine longest length of increasing subsequences in permutations~\cite{kiwi4}  and various properties of the largest size of a planar matching in random bipartite graphs have been studied in~\cite{joh}~\cite{kiwi}~\cite{kiwi2}~\cite{kiwi3}.

For example, the paper~\cite{joh} studied weighted planar matchings of random bipartite graphs with geometrically distributed edge weights and obtained expressions for the asymptotic expected weight. The paper~\cite{kiwi} considered dependent planar matchings (i.e., where the edge states are not independent of one another) with additional degree constraints and obtained convergence of the maximum size, appropriately scaled, and in~\cite{kiwi2}, the analysis is extended to general planar subgraphs of not necessarily planar graphs.

In this paper, we study planar matchings from a probabilistic perspective. Our motivation comes from analyzing redundancy in datasets via the Shapley value allocation~\cite{shap} as an alternative to the graph based approach described in~\cite{ganesan2}. In such a setup, edge weights correspond to the redundancy contribution of a data point in a given ordering and planar matchings correspond to subsets of dissimilar data points. It is of interest to estimate the minimum redundancy of such a data subset.

In what follows, we study minimum weight planar matchings in bipartite graphs equipped with positive random edge weights. We use segmentation techniques to estimate the minimum weight of planar matchings with a \emph{given} number of edges in Theorem~\ref{thm_weight_ax}. Finally, in Theorem~\ref{thm_weight_two}, we describe a saturation property of the minimum weight that occurs when the constraint on the number of edges in the planar matching becomes comparable to the number of vertices.

The paper is organized as follows: In Section~\ref{sec_max_size}, we obtain bounds for the maximum size of planar matchings with edge length constraints and in Section~\ref{sec_min_wt}, we describe and study minimum weight planar matchings and obtain deviation and convergence estimates for the minimum weight.

\setcounter{equation}{0}
\renewcommand\theequation{\thesection.\arabic{equation}}
\section{Length Constrained Planar Matchings} \label{sec_max_size}
For~\(i \geq 1,\) let~\(u_i = (i,1)_2\) and~\(v_i = (i,0)_2\) be points in~\(\mathbb{R}^2.\) We use the subscript~\(2\) to differentiate from the two-tuple notations for edges introduced later. For each~\(i \neq j,\) join the vertices~\(u_i\) and~\(v_j\) by an edge (denoted as~\((u_i,v_j)\)) to obtain an infinite bipartite graph~\(K_{tot}.\)

For integer~\(n \geq 1,\) let~\(K_{n,n}\) be the complete bipartite graph containing~\(n\) top vertices~\(X_{n,n}=\{u_1,\ldots,u_n\}\) and~\(n\) bottom vertices~\(Y_{n,n} = \{v_1,\ldots,v_n\}\) and let~\(E(K_{n,n})\) be the edge set of~\(K_{n,n}.\) An edge~\(e  \in E(K_{n,n})\) if and only if~\(e\) has one end-vertex~\(u_i \in X_{n,n}\) and the other end-vertex~\(v_j \in Y_{n,n}.\) We denote the edge~\(e\) as~\(e = (u_i,v_j)\) and  define the \emph{length} of~\(e\) to be~\(|i-j|.\) This is purely for convenience and our results and an analogous analysis as below holds if we use the usual Euclidean length instead.

\begin{figure}[tbp]
\centering
\includegraphics[width=3in, trim= 20 450 20 120, clip=true]{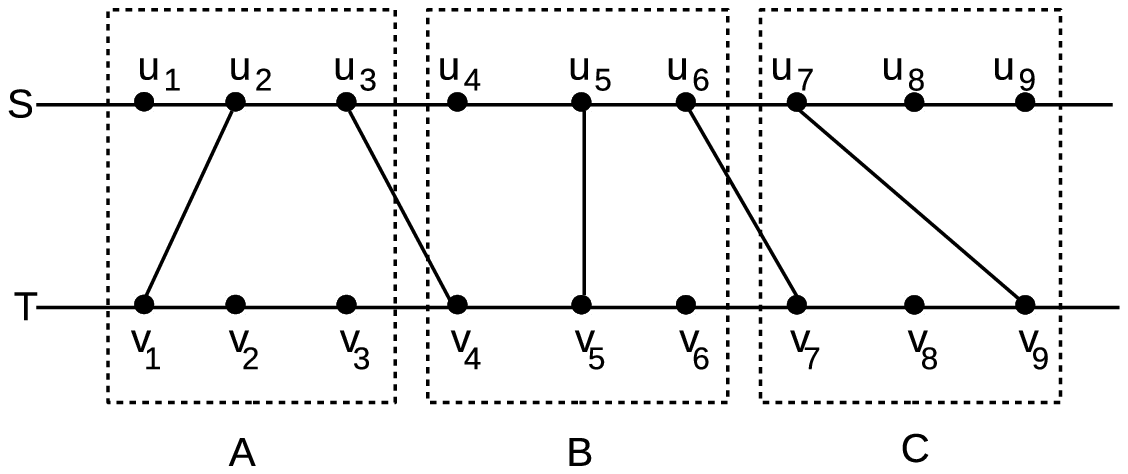}
\caption{Illustration of a planar matching.}
\label{fig_match_ax}
\end{figure}

A \emph{matching} of size~\(t\) in~\(K_{n,n}\) is a set of vertex disjoint edges~\({\cal M} = \{e_1,\ldots,e_t\}.\) Suppose~\(e_l = (u_{i_l},v_{j_l})\) for~\(1 \leq l \leq t.\)  We say that~\({\cal M}\) is a \emph{planar matching} if~\(i_1<i_2< \ldots < i_t\) and~\(j_1< j_2 < \ldots < j_t.\) In other words, no two edges in~\({\cal M}\) intersect each other. The size of~\({\cal M}\) is defined to be the number of edges~\(t.\)
\begin{definition}\label{plan_mat_def}
For~\(L \geq 0,\) we say that~\({\cal M}\) is an~\(L-\)constrained planar matching if the length of each edge in~\({\cal M}\) is at most~\(L;\) i.e.,~\(|i_k-j_k| \leq L\) for each~\(1 \leq k \leq t.\)
\end{definition}
In Figure~\ref{fig_match_ax}, we illustrate the above definition with an example of a planar matching~\({\cal M}_0\) containing five edges~\((u_2,v_1),(u_3,v_4),(u_5,v_5), (u_6,v_7)\) and~\((u_7,v_9).\) The edge~\((u_7,v_9)\) has the largest length equal to~\(2\) and so~\({\cal W}_0\) is a~\(2-\)constrained planar matching. 

Let~\(\{Y(f)\}_{f \in E(K_{n,n})}\) be independent and identically distributed (i.i.d.) Bernoulli random variables indexed by the edge set of~\(K_{n,n}\) and having distribution
\begin{equation}\label{x_dist}
\mathbb{P}(Y(f)= 1) = p(f) =1-\mathbb{P}(Y(f) = 0)
\end{equation}
for each edge~\(f\) and some~\(0 \leq p(f)\leq 1.\) We refer to~\(Y(f)\) as the \emph{state} of the edge~\(f\) and denote~\(G\) to be the random subgraph of~\(K_{n,n}\) formed by the set of all edges satisfying~\(Y(f) =1.\) In general, the edge probability~\(p(f)\) need not be the same for all edges and in that case, we say that the random graph~\(G\) obtained is \emph{inhomogenous}. If~\(p(f) = p\) for all~\(f\) and some~\( 0 \leq p \leq 1\) then we say that~\(G\) is a \emph{homogenous} random graph.

Letting~\(U_{n} = U_n(L)\) be the maximum size of an~\(L-\)constrained planar matching in~\(G,\) we are interested in studying the growth of~\(U_n\) as a function of~\(n,L\) and the average edge probability defined as follows: we let~\({\cal T}_L\) be the set of all edges of~\(K_{n,n}\) of length at most~\(L\) and  define
\begin{equation}\label{t_cond}
p_{av} = p_{av}(L) := \frac{1}{\zeta_L} \sum_{f \in {\cal T}_L} p(f),
\end{equation}
to be the corresponding averaged edge probability, where~\(\zeta_L\) is the size of~\({\cal T}_L.\) We have the following deviation and variance estimates for~\(U_n.\) Throughout, constants do not depend on~\(n.\)
\begin{theorem}\label{prop_one} We have
\begin{equation}\label{var_est_one}
var(U_n) \leq 2 \mathbb{E}U_n
\end{equation}
and
\begin{equation}\label{mn_inhom}
\mathbb{P}\left(\frac{np_{av}}{28L} \leq M_n \leq \frac{9np_{av}L}{2}\right) \geq 1-2\exp\left(-\frac{np_{av}L}{16}\right).
\end{equation}
\end{theorem}
Thus if~\(\frac{np_{av}}{L} \longrightarrow \infty,\) then from~(\ref{mn_inhom}) we get that~\[\mathbb{E}U_n \geq \frac{np_{av}}{28L} \left(1-e^{-\frac{3np_{av}L}{32}}\right) \longrightarrow \infty\] as~\(n \rightarrow \infty.\) From the variance estimate~(\ref{var_est_one}), we then obtain~\(L^2-\)convergence:~\[\mathbb{E}\left(\frac{U_n}{\mathbb{E}U_n}-1\right)^2  \leq \frac{2}{\mathbb{E}U_n} \longrightarrow 0\]  as~\(n \rightarrow \infty.\)

Throughout, we use the following standard deviation estimate regarding sums of independent Bernoulli random variables. Let~\(\{X_j\}_{1 \leq j \leq r}\) be independent Bernoulli random variables with~\(\mathbb{P}(X_j = 1) = 1-\mathbb{P}(X_j = 0) > 0.\) If~\(T_r := \sum_{j=1}^{r} X_j,\theta_r := \mathbb{E}T_r\) and~\(0 < \gamma \leq \frac{1}{2},\) then
\begin{equation}\label{conc_est_f}
\mathbb{P}\left(\left|T_r - \theta_r\right| \geq \theta_r \gamma \right) \leq 2\exp\left(-\frac{\gamma^2}{4}\theta_r\right)
\end{equation}
for all \(r \geq 1.\) For a proof of~(\ref{conc_est_f}), we refer to Corollary A.1.14, pp. 312 of~\cite{alon}.

\emph{Proof of Theorem~\ref{prop_one}}: The variance estimate~(\ref{var_est_one}) is obtained using the standard martingale difference method (see for example, Theorem~\(2.2\)  of~\cite{ganesan}). For completeness, we give a quick proof here. Let~\(\{e_j\}_{1 \leq j \leq n^2}\) be the set of edges of~\(K_{n,n}\) and let~\(U_n^{(j)}\) be the maximum size of a planar matching in the graph~\(G^{(j)}\) obtained by replacing the state~\(Y(e_j)\) of the edge~\(e_j\) in the graph~\(G\) by an independent copy~\(Y^{(c)}(e_j),\) while retaining the state of every other edge. From the Efron-Stein inequality (see Section~\(2,\) Eq.~\((2.1)\) of~\cite{steele}) we then get
\begin{equation}\label{var_un}
var(U_n) \leq \sum_{j=1}^{n^2} \mathbb{E}\left(U_n-U_n^{(j)}\right)^2.
\end{equation}

Changing the state of the edge~\(e_j\) changes the maximum size of an~\(L-\)constrained planar matching by at most one and moreover a change occurs only if~\(e_j\) belongs to the maximum size~\(L-\)constrained planar matching  in either~\(G\) or~\(G^{(j)}.\) In other words, if~\({\cal W}\) and~\({\cal W}^{(j)}\) denote the maximum size~\(L-\)constrained planar matchings in~\(G\) and~\(G^{(j)},\) respectively, then
\begin{align}
|U_n-U_n^{(j)}| &\leq \ind\left( \left\{e_j \in {\cal W} \right\} \bigcup \left\{e_j \in {\cal W}^{(j)}\right\}\right) \nonumber\\
&\leq \ind\left(e_j \in {\cal W}\right) + \ind\left(e_j \in {\cal W}^{(j)}\right). \label{un_diff}
\end{align}

Squaring~(\ref{un_diff}), taking expectations and using~\((a+b)^2 \leq 2(a^2+b^2),\) we get that
\begin{equation}
\mathbb{E}\left(U_n - U_n^{(j)}\right)^2 \leq 2\mathbb{P}\left(e_j \in {\cal W}\right) + 2\mathbb{P}\left(e_j \in {\cal W}^{(j)}\right) = 4\mathbb{P}\left(e_j \in {\cal W}\right), \label{un_diff_two}
\end{equation}
since~\({\cal W}\) and~\({\cal W}^{(j)}\) are identically distributed. Plugging~(\ref{un_diff_two}) into~(\ref{var_un}), we get
\[var(U_n) \leq 4\sum_{j=1}^{n^2} \mathbb{P}\left(e_j \in {\cal W}\right) = 4\mathbb{E}U_n,\] since~\({\cal W}\) has size~\(U_n,\) by definition. This completes the proof of~(\ref{var_est_one}).

Next, to prove~(\ref{mn_inhom}), let~\({\cal S}_L = \{f_1,\ldots,f_Q\} \subseteq {\cal T}_L\) be the set of all edges of the random graph~\(G,\) each of length at most~\(L.\) We consider the graph~\({\cal H}_L\) obtained as follows. The edge~\(f_i\) is represented by a vertex~\(z_i\) in the graph~\({\cal H}_L.\) Two vertices~\(z_i\) and~\(z_j\) are joined by an edge in~\({\cal H}_L\) if and only if the edges~\(f_i\) and~\(f_j\) intersect in~\(K_{n,n}.\)

A stable set in~\({\cal H}_L\) is a set of vertices no two of which are adjacent to each other and by construction, any stable set in~\({\cal H}_L\) corresponds to an~\(L-\)constrained planar matching in~\(G.\) If~\(m\) denotes the number of edges in~\({\cal H}_L,\) then the average vertex degree in~\({\cal H}_L\) is~\(d_{av} := \frac{2m}{Q}\) and from Theorem~\(3.2.1,\) pp.~\(27\) of~\cite{alon} we know that~\({\cal H}_L\) contains a stable set of size at least~\[\frac{Q}{2\max(1,d_{av})} \geq \frac{Q}{2(d_{av}+1)} = \frac{Q^2}{2(2m+Q)}\] and so
\begin{equation}\label{mn_low}
\frac{Q^2}{2(2m+Q)} \leq U_n \leq Q.
\end{equation}

In what follows, we estimate the size~\(Q\) of~\({\cal S}_L\) and~\(m\) in that order. From the lower bound in condition~(\ref{t_cond}),
we see that~\(\mathbb{E}Q = p_{av} \cdot \zeta_L,\) where we recall from~(\ref{t_cond}) that~\(\zeta_L\) is the number of edges of~\(K_{n,n}\) of length at most~\(L.\)  Therefore using the deviation estimate~(\ref{conc_est_f}) with~\(\gamma = \frac{1}{2},\) we get that
\begin{equation}\label{q_est}
\mathbb{P}\left(\frac{p_{av}\zeta_L}{2}  \leq Q \leq \frac{3p_{av}\zeta_L}{2} \right) \geq 1-2\exp\left(-\frac{p_{av}\zeta_L}{16}\right).
\end{equation}
We provide a quick estimate of~\(\zeta_L\) as follows. If an edge~\(e \in E(K_{n,n})\) has end-vertices~\(u_i \in X\) and~\(v_j \in Y\) and the length of~\(e\) is at most~\(l,\) then necessarily~\(|i-j| \leq l.\) Therefore the number of edges of length exactly equal to~\(l \neq 0\) is~\(2(n-l)\) and moreover, the number of edges of length zero, is exactly~\(n.\) Consequently
\begin{equation}\label{tau_L_est}
\zeta_L = n + \sum_{1 \leq i \leq L} 2(n-i) =  nL + (n-L)(L+1)
\end{equation}
and so~\(nL \leq \zeta_L \leq 3nL.\) Plugging this into~(\ref{q_est}) we therefore get that
\begin{equation}\label{q_est2}
\mathbb{P}\left( \frac{np_{av}L}{2}  \leq Q \leq \frac{9np_{av}L}{2}\right) \geq 1-2\exp\left(-\frac{np_{av}L}{16}\right).
\end{equation}

\begin{figure}[tbp]
\centering
\includegraphics[width=3in, trim= 20 490 20 80, clip=true]{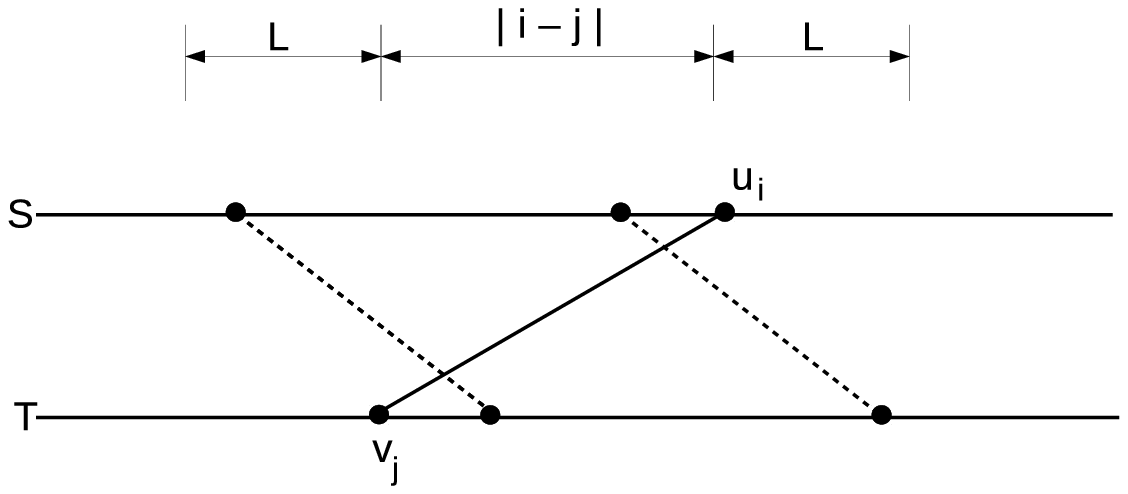}
\caption{Determining the range of edges that intersect the edge~\(e\) with end-vertices~\(u_i\) and~\(v_j.\)}
\label{fig_match_ax5}
\end{figure}

To estimate~\(m,\) we use the identity that the number of edges in any graph is twice the sum of all vertex degrees. Any vertex in~\(K_{n,n}\) is adjacent to at most~\(2L\) edges of length at most~\(L.\) Therefore if an edge~\(e \in G\) has length~\(l \leq L,\) there are at most~\((l+2L) \cdot 2L \leq 6L^2\) edges in~\(K_{n,n}\) that intersect~\(e.\) To see this is true suppose~\(e\) has end-vertices~\(u_i \in X\) and~\(v_j \in Y\) with~\(i \geq j.\) Any edge of length at most~\(L\) that intersects~\(e\) must have end-vertices of the form~\(u_x\) and~\(v_y,\) where~\( i-l-L\leq x \leq i+L\) and~\( j-L \leq y \leq j+l+L.\) This is illustrated in Figure~\ref{fig_match_ax5} where we show two possible edges (shown by dotted lines) that intersect the edge~\(e.\)  Therefore, the sum of vertex degrees in the graph~\({\cal H}_L\) is at most \(6QL^2\) and this implies that~\(m \leq 3QL^2.\)

Consequently, if~\(Q\) satisfies the bounds in~(\ref{q_est2}), then we get that
\begin{equation}\label{q_chain}
\frac{Q^2}{2(Q+2m)} \geq \frac{Q^2}{2(Q+6QL^2)} \geq \frac{Q^2}{14QL^2}  = \frac{Q}{14L^2} \geq \frac{np_{av}}{28L}
\end{equation}
where the second inequality in~(\ref{q_chain}) follows from the fact that~\(L \geq 1.\) From~(\ref{mn_low}),~(\ref{q_est2}) and~(\ref{q_chain}), we then get~(\ref{mn_inhom}) and this completes the proof of the Theorem.~\(\qed\)

\setcounter{equation}{0}
\renewcommand\theequation{\thesection.\arabic{equation}}
\section{Minimum Weight Planar Matchings} \label{sec_min_wt}
Recall the infinite bipartite graph~\(K_{tot}\) defined in the previous section and let~\(w(f), f \in K_{tot}\) be independent (but not necessarily identically distributed) random variables, indexed by the edge set of~\(K_{tot}\) and satisfying~\(w(f) \geq 0\) a.s.\ for each edge~\(f.\) We define~\(w(f)\) to be the \emph{weight} of the edge~\(f.\)

As before, let~\(K_{n,n} \subset K_{tot}\) be the complete bipartite graph with top vertex set~\(X_{n,n}  =\{u_1,\ldots,u_n\}\) and bottom vertex set~\(Y_{n,n} = \{v_1,\ldots,v_n\}.\) For a planar matching~\({\cal P}\) we define~\(W({\cal P}) := \sum_{f \in {\cal P}} w(f)\) to be the sum of the weights of edges in~\({\cal P}.\)
\begin{definition} For~\(1\leq \tau \leq n,\) we define~\(M_n = M_n(\tau)\) to be the minimum weight of a planar matching of~\(K_{n,n},\) containing at least~\(\tau\) edges.
\end{definition}
We are primarily interested in studying how~\(M_n(\tau)\) behaves in terms of the size constraint~\(\tau\) and the edge weight distribution parameters.





To begin with, we assume that there are non-increasing functions\\\(0 \leq F_{low}(x) \leq F_{up}(x) \leq 1\) such
that
\begin{equation}\label{g_up_def_ax}
F_{low}(x) \leq \inf_{f} \mathbb{P}(w(f) \leq x)  \leq \sup_{f} \mathbb{P}(w(f) \leq x) \leq F_{up}(x)
\end{equation}
and satisfying~\(F_{up}(0) = 0\) and~\(0 < F_{low}(x) \leq 1\) for each~\(0 < x \leq 1.\) In words,~\(F_{low}\) and~\(F_{up}\) respectively are lower and upper bounds on the cumulative distribution function (cdf) of the edge weights and in the special case where the weights are independent and identically distributed (i.i.d.) with common cdf~\(F,\) we have~\(F_{low} = F_{up} = F.\)

With the above notation, the terms
\begin{equation}\label{an_bn_def_ax}
s_n := \inf\left\{x >0 : F_{up}(x) \geq \left(\frac{\tau}{64n}\right)^2 \right\}\text{ and }t_n := \sup\left\{x >0 : F_{low}(x) \leq  \left(\frac{8\tau}{n}\right)^2\right\}
\end{equation}
are essentially quantiles of the order of~\(\frac{\tau^2}{n^2}\) and we have the following result regarding the expectation and deviation behaviour of~\(M_n\) as a function of the edge constraint~\(\tau, s_n\) and~\(t_n.\) Throughout, constants do not depend on~\(n\) or~\(\tau.\)
\begin{theorem}\label{thm_weight_ax} Suppose~\(\sup_{f} \mathbb{E}w^2(f) < \infty\)  and let~\((\log{n})^6 \leq \tau \leq n.\)\\
\((a)\) There are  constants~\(\alpha_1,\alpha_2 >0\) such that
\begin{equation}\label{dev_wn_low_ax}
\mathbb{P}(M_n \leq \alpha_1 \tau t_n) \geq 1-e^{-\alpha_2 \tau} \text{ and }\mathbb{E}M_n \leq \alpha_1 \tau t_n + \alpha_1 n^2e^{-\alpha_2 \tau}.
\end{equation}
\((b)\) There are constants~\(\beta_1,\beta_2 > 0\) such that
\begin{equation}\label{dev_wn_up_ax}
\mathbb{P}\left(M_n \geq \beta_1  \tau s_n\right) \geq 1-e^{-\beta_2(\log{n})^2} \text{ and } \mathbb{E}M_n \geq \beta_1 \tau s_n.
\end{equation}
\((c)\) There is a constant~\(\gamma > 0\) such that
\begin{equation}\label{var_bound_ax}
var(M_n)  \leq \gamma c_n,
\end{equation}
where~\(c_n := \min\left(n,\frac{D \tau t_n}{s_n}\right),\) the constant~\(D = \max\left(1,\frac{2\alpha_1}{\beta_1}\right)\) and~\(\alpha_1,\beta_1\) are as in~(\ref{dev_wn_low_ax}) and~(\ref{dev_wn_up_ax}), respectively.
\end{theorem}
The above result obtains deviation, expectation and variance bounds for the minimum weight of constrained planar matchings, assuming that the edge weights have bounded second moment. These could be used to investigate~\(L^2-\)convergence of~\(M_n\) as illustrated below.

For example, suppose the weights are independent and identically distributed (i.i.d.) with~\(F_{low}(x) = F_{up}(x) = x^{\alpha}\) for~\(0 \leq x \leq 1\) and some constant~\(\alpha > 0.\) Thus the edge weights are bounded a.s.\ and so trivially satisfy the bounded second moment condition. From the definitions of~\(s_n\) and~\(t_n\) in~(\ref{an_bn_def_ax}), we see that
\[s_n = \left(\frac{\tau}{64n}\right)^{2/\alpha} \text{ and } t_n = \left(\frac{8\tau}{n}\right)^{2/\alpha}\] so that,~\(s_n\) and~\(t_n\) are both of the same order~\(\left(\frac{\tau}{n}\right)^{2/\alpha}.\) From the deviation bounds~(\ref{dev_wn_low_ax}),~(\ref{dev_wn_up_ax}) and the union bound, we then get for~\((\log{n})^6 \leq \tau \leq n\) that
\begin{align}
&\mathbb{P}\left( C_1 \tau \left(\frac{\tau}{n}\right)^{2/\alpha}  \leq M_n \leq C_2 \tau \left(\frac{\tau}{n}\right)^{2/\alpha} \right) \nonumber\\
&\;\;\;\;\;\;\;\;\geq\;\;\;1-\exp\left(-C_3 (\log{n})^2\right) - \exp\left(-C_3 (\log{n})^6\right) \nonumber\\
&\;\;\;\;\;\;\;\;\geq\;\;\;1-2\exp\left(-C_3 (\log{n})^2\right) \nonumber
\end{align}
for some constants~\(C_1,C_2,C_3> 0.\)

Also if~\((\log{n})^6 \leq \tau \leq n,\) then the term~\(\tau t_n = \tau \left(\frac{\tau}{n}\right)^{2/\alpha}\) is much larger than~\(n^2e^{-C\tau}\) for any constant~\(C > 0\) and all~\(n\) large and so from the expectation bounds in~(\ref{dev_wn_low_ax}) and~(\ref{dev_wn_up_ax}), we get that
\begin{equation}\label{exp_mn_bounds}
D_1 \tau \left(\frac{\tau}{n}\right)^{2/\alpha} \leq \mathbb{E}M_n \leq D_2 \tau \left(\frac{\tau}{n}\right)^{2/\alpha}
\end{equation}
for some constants~\(D_1,D_2 > 0.\) Finally again using the fact that~\(s_n\) and~\(t_n\) are of the same order~\(\left(\frac{\tau}{n}\right)^{2/\alpha},\) we see that the term~\(c_n\) in the variance estimate~(\ref{var_bound_ax}) is bounded above as~\(c_n \leq D_3 \min(n,\tau) \leq D_3 \tau\) for some constant~\(D_3 > 0.\)  Thus from~(\ref{var_bound_ax}), we see that~\(var(M_n) \leq D_4 \tau\) for some constant~\(D_4 > 0\) and combining with the lower bound for~\(\mathbb{E}M_n\) in~(\ref{exp_mn_bounds}), we get that
\[\mathbb{E}\left(\frac{M_n}{\mathbb{E}M_n}-1\right)^2 \leq \frac{D_4 \tau}{D_1^2 \tau^2 \left(\tau/n\right)^{4/\alpha}} = D_5 \frac{n^{4/\alpha}}{\tau^{4/\alpha+1}} \longrightarrow 0\] for some constant~\(D_5 > 0,\) provided~\(n^{\theta} \leq \tau \leq n\) for some~\(\theta > \frac{4}{4+\alpha},\) strictly.

\underline{\emph{Remark}}: For~\(\tau = o(n),\) that is if~\(\frac{\tau}{n} \longrightarrow 0\) then from~(\ref{exp_mn_bounds}), we see that~\(\mathbb{E}M_n(\tau)\) exhibits different behaviour for different values of~\(\alpha\) and so is dependent on the decay of the cumulative distribution function (cdf) of the edge weight distributions near the origin. However, if~\(\tau\) is of the order of~\(n,\) the minimum weight~\(M_n\) is also of the order of~\(n\) with high probability, i.e., with probability converging to one as~\(n \rightarrow \infty\) irrespective of the cdf.  In other words,~\(M_n(\tau)\) \emph{saturates} when~\(\tau\) changes from being sub-linear to growing linearly with~\(n.\) We study this aspect in our next result at the end of this section.

\underline{\emph{Remark 2}}: In this paper, we have studied variance and~\(L^2-\)convergence of the minimum weight~\(M_n.\) For future, it might be interesting to investigate Central Limit Theorems for~\(M_n.\) New methods might be needed to obtain variance lower bounds and thereby invoke Lindeberg-Feller type conditions for the same.




\emph{Proof of Theorem~\ref{thm_weight_ax}\((a)\)}: We use a segmentation approach as in the proof of Theorem~\(1,\) Section~\(6.1\) of~\cite{kiwi4}, which studies the related but dual problem of determining the length of the \emph{longest} common subsequence in large alphabets. Assuming for simplicity that~\(\frac{n}{4\tau}\) is an integer, we divide the edge set of~\(K_{n,n}\) into~\(4\tau\) disjoint segments each containing~\(\left(\frac{n}{4\tau}\right)^2\) edges. Formally, let~\({\cal A}_i, 1 \leq i \leq 4\tau\) be the set of all top vertices~\(u_j\) with~\((i-1)\frac{n}{4\tau} + 1 \leq j \leq \frac{i n}{4\tau}\) and let~\({\cal B}_i, 1 \leq i \leq 4\tau\) be the set of all bottom vertices~\(v_j\) with~\((i-1)\frac{n}{4\tau} + 1 \leq j \leq \frac{i n}{4\tau}.\) For example in Figure~\ref{fig_match_ax},~\(n=9\) and we have split~\(K_{9,9}\) into three segments~\(A,B\) and~\(C\) shown as dotted rectangles.

An edge~\(f\) is said to belong to the~\(i^{th}\) segment if~\(f\) has one endvertex in~\({\cal A}_i\) and the other endvertex in~\({\cal B}_i.\) Say that~\(f\) is \emph{good} if the weight~\(w(f) \leq 2t_n\) and \emph{bad} otherwise, where~\(t_n\) is as defined in~(\ref{an_bn_def_ax}). The probability that an edge~\(f\) in the~\(i^{th}\) segment is bad, is at most~\(1-\left(\frac{8\tau}{n}\right)^2\) and so the probability that each edge in the~\(i^{th}\) segment is bad is bounded above by~\[\left(1- \left(\frac{8\tau}{n}\right)^2\right)^{\left(\frac{n}{4\tau}\right)^2} \leq e^{-4} \leq \frac{1}{32}.\] Therefore if~\(F_i\) denotes the event that at least one edge in the~\(i^{th}\) segment is good, then~\(\mathbb{P}(F_i) \geq \frac{31}{32}.\)

Letting~\(X := \sum_{i=1}^{4\tau} \ind(F_i)\) be the total number of segments containing at least one good edge, we  get that~\[\mathbb{E}X \geq \frac{31\tau}{8} \geq \tau.\] For~\(i \neq j,\) the events~\(F_i\) and~\(F_j\) are independent and so using the deviation estimate~(\ref{conc_est_f}) with~\(\gamma = \frac{1}{2},\) we see that
\begin{equation}\label{x_est}
\mathbb{P}\left( X \geq \frac{31\tau}{16} \right) \geq 1-e^{-2C_1\tau}
\end{equation}
for some constant~\(C_1 > 0.\) If~\(X \geq \frac{31\tau}{16}> \tau, \) then there at least~\(\tau\) segments, each containing a good edge. These good edges form a planar matching containing at least~\(\tau\) edges and having weight at most~\(2\tau t_n.\) This obtains the deviation upper bound for~\(M_n\) in~(\ref{dev_wn_low_ax}).

For estimating~\(\mathbb{E}M_n,\) we use the upper bound~\(M_n \leq \sum_{f \in E(K_{n,n})} w(f) \) whenever the event~\(E_{bad} := \left\{X  < \frac{31\tau}{16}\right\}\) occurs, to get that
\begin{equation}\label{waldo}
\mathbb{E}M_n \leq 2\tau t_n + \mathbb{E}\sum_{f \in E(K_{n,n})}w(f) \ind( E_{bad}).
\end{equation}
Using the Cauchy-Schwarz inequality, we get for any edge~\(f\) that
\begin{equation}\label{daldo}
\mathbb{E}w(f) \ind(E_{bad}) \leq \left(\mathbb{E}w^2(f)\right)^{1/2} \mathbb{P}^{1/2}\left(E_{bad}\right) \leq C_2e^{-C_1\tau},
\end{equation}
for some constant~\(C_2 > 0\) where the final estimate in~(\ref{daldo}) follows from~(\ref{x_est}) and the fact that the edge weights have uniformly bounded second moments (see Theorem statement).

Plugging~(\ref{daldo}) into~(\ref{waldo}) and using the fact that there are~\(n^2\) edges in~\(K_{n,n},\) we obtain the upper expectation bound in~(\ref{dev_wn_low_ax}). This completes the proof of Theorem~\ref{thm_weight_ax}\((a).\)~\(\qed\)

\emph{Proof of Theorem~\ref{thm_weight_ax}\((b)\)}: We use the following upper bound for maximal Bernoulli planar matchings: Suppose each edge of the complete bipartite graph~\(K_{n,n}\) is present with probability at most~\(p_{up},\) independent of the other edges. If~\(R_n\) is the maximum \emph{size} of a planar matching in~\(K_{n,n},\) then
\begin{equation}\label{mn_est}
\mathbb{P}\left(R_n \geq 4en \sqrt{p_{up}}\right) \leq \frac{1}{2^{4en\sqrt{p_{up}}}} \text{ and } \mathbb{E}R_n \leq 16en \sqrt{p_{up}}.
\end{equation}
For completeness, we provide a small proof  of~(\ref{mn_est}) in the Appendix and for sharper upper and lower bounds for~\(R_n,\) we refer to Theorem~\(11\) of~\cite{kiwi4}.

Let~\(s_n\) be as in~(\ref{an_bn_def_ax}) and colour an edge~\(f \in E(K_{n,n})\) to be red if its weight~\(w(f) \leq \frac{s_n}{2}\) and blue otherwise. By  definition of~\(s_n,\) we see that~\(f\) has a red colour with probability  at most~\(p_{up} = \left(\frac{\tau}{64n}\right)^2.\) Therefore from the upper expectation bound in~(\ref{mn_est}) above, we get that the expected maximum size of a planar matching formed by red edges is at most~\(\frac{e\tau}{4}.\)

Let~\({\cal M}\) be a minimum weight planar matching containing  at least~\(\tau\) edges. By definition, the weight of~\({\cal M}\) is~\(M_n(\tau)\) and we let~\(M_{red}\) be the total weight of the red edges in~\({\cal M}.\) The red edges in~\({\cal M}\) themselves form a planar  matching and from the above paragraph, we see that the expected number of red edges in~\({\cal M}\) is at most~\(\frac{e\tau}{4}.\) Each red edge has weight at most~\(\frac{s_n}{2}\) and so~\[\mathbb{E}M_{red} \leq \frac{e \tau}{4} \cdot \frac{s_n}{2}.\] The rest of the edges in~\({\cal M}\) have weight at least~\(\frac{s_n}{2}\) and so~\[\mathbb{E}M_n(\tau) \geq \left(\tau - \frac{e\tau}{4}\right) \frac{s_n}{2}.\] This proves the desired expectation lower bound for~\(M_n.\)

To obtain the deviation lower bound for~\(M_n,\) we again estimate the total weight of red edges as before, with an additional step to take care of ``low probability" edges. Formally, mark an edge~\(f \in E(K_{n,n})\) with a cross symbol if
\begin{equation}\label{mark_def}
\mathbb{P}\left(w(f) \leq \frac{s_n}{2}\right) \geq \frac{(\log{n})^4}{n^2}
\end{equation}
and as before, let~\({\cal M}\) be a minimum weight planar matching containing  at least~\(\tau\) edges. We recall that~\(M_{red}\) is the total weight of red edges in~\({\cal M}\) and upper bound~\(M_{red}\) as follows. If~\(M_{red}(1) \leq M_{red}\) is the total weight of \emph{marked} red edges within~\({\cal M}\) and~\(M_{red}(2)\) is the total weight of all \emph{unmarked} red edges in~\(K_{n,n},\) then
\begin{equation}\label{jalpa}
M_{red} \leq M_{red}(1) + M_{red}(2).
\end{equation}
We estimate~\(M_{red}(1)\) and~\(M_{red}(2)\) in that order  below and then combine the respective estimates together at the end.

If~\(f\) is a red edge then~\(w(f) \leq \frac{s_n}{2}\) and so if~\[p_{up} := \max_{f \text{ marked}} \mathbb{P}\left(w(f) \leq \frac{s_n}{2}\right)\] is the maximum probability of occurrence of a marked edge, then from  the definition of~\(s_n\) in~(\ref{an_bn_def_ax}) we  see that
\[\frac{(\log{n})^4}{n^2} \leq p_{up} \leq \left(\frac{\tau}{64n}\right)^2,\] provided there is at least one marked edge. Denoting~\(N_{red}(1)\) to be the \emph{number} of marked red edges in~\({\cal M}\) and using the upper deviation bound in~(\ref{mn_est}) stated in the first paragraph of this proof, we then get that
\begin{equation}\label{num_red_est}
\mathbb{P}\left(N_{red}(1) \geq \frac{e\tau}{4}\right) \leq e^{-2C(\log{n})^2}
\end{equation}
for some constant~\(C > 0.\) By definition, the total weight of marked red edges is~\[M_{red}(1) \leq N_{red}(1) \frac{s_n}{2}\] and so
\begin{equation}\label{m_red_one_est}
\mathbb{P}\left(M_{red}(1) \geq \frac{e\tau s_n}{8}\right) \leq e^{-2C(\log{n})^2}.
\end{equation}

Next, for an unmarked edge~\(f\) let~\[X_f = \ind(f \text{ is red}) = \ind\left(w(f) \leq \frac{s_n}{2}\right)\] be the indicator function that~\(f\) is red and let~\(q_f = \mathbb{P}(X_f=1).\) From the definition of marked edge in~(\ref{mark_def}), we see that~\(q_f \leq \frac{(\log{n})^4}{n^2}\) and so
\[\mathbb{E}e^{X_f}  = (1-q_f + eq_f) \leq \exp(q_f(e-1)) \leq \exp\left((e-1)\frac{(\log{n})^4}{n^2}\right).\] If~\(N_{red}(2)\) is the \emph{total} number of unmarked red edges in~\(K_{n,n},\) then from the Chernoff bound we get that
\begin{align}\label{third_term_ax}
\mathbb{P}\left(N_{red}(2) \geq 2(\log{n})^5\right) &\leq e^{-2(\log{n})^5} \prod_{f \in E(K_{n,n})} \mathbb{E}e^{X_f} \nonumber\\
&\leq e^{-2(\log{n})^5}  \cdot  \exp\left((e-1)(\log{n})^4\right) \nonumber\\
&\leq e^{-(\log{n})^5}
\end{align}
for all~\(n\) large, where the second inequality is true since there are~\(n^2\) edges in~\(K_{n,n}.\) If~\(N_{red}(2) \leq 2(\log{n})^5\) then the total weight~\(M_{red}(2)\) of all unmarked red edges is at most~\(2(\log{n})^5 \cdot \frac{s_n}{2} = s_n (\log{n})^5\) and so from~(\ref{third_term_ax}) we see that
\begin{equation}\label{m_red_two_est}
\mathbb{P}\left(M_{red}(2) \geq (\log{n})^5 s_n\right) \leq e^{-(\log{n})^5}.
\end{equation}

Combining~(\ref{m_red_one_est}) and~(\ref{m_red_two_est}) and using~(\ref{jalpa}), we get that
\[\mathbb{P}\left(M_{red} \geq \left(\frac{e\tau}{4} + 2(\log{n})^5\right)\frac{s_n}{2}\right) \leq e^{-2C(\log{n})^2} + e^{-(\log{n})^5} \leq e^{-C(\log{n})^2}.\] Any edge that is not red has weight at least~\(\frac{s_n}{2}\) and so the total weight of all the edges in the minimum weight planar matching~\({\cal M}\) containing at least~\(\tau\) edges, is at least~\[\left(\tau - \frac{e\tau}{4} - 2(\log{n})^5\right) \frac{s_n}{2}\] with probability at least~\(1-e^{-C(\log{n})^2}.\) Since~\(\tau \geq (\log{n})^6,\) this obtains the lower deviation bound for~\(M_n\) in~(\ref{dev_wn_low_ax}) and therefore completes the proof of Theorem~\ref{thm_weight_ax}\((b).\)~\(\qed\)


\emph{Proof of Theorem~\ref{thm_weight_ax}\((c)\)}: We use the martingale difference method to obtain the variance estimate. Let~\(e_j, 1 \leq j \leq n^2\) be the set of edges of  the complete graph~\(K_{n,n}\) and let~\(M_n^{(j)}\) be the minimum weight of a planar matching containing at least~\(\tau\) edges, when the weight~\(w(e_j)\) of the edge~\(e_j\) is replaced by an independent copy~\(w^{(c)}(e_j),\) leaving all other edge weights unchanged. From Efron-Stein inequality (see Section~\(2,\) Eq.~\((2.1)\) of~\cite{steele}), we then get that
\begin{equation}\label{efron_ax}
var(M_n) \leq \sum_{j=1}^{n^2} \mathbb{E}\left(M_n-M_n^{(j)}\right)^2.
\end{equation}

To estimate the difference~\(|M_n-M_n^{(j)}|,\) let~\({\cal M}\) be the minimum weight planar matching with edge weights~\(\{w(e_j)\}_{1 \leq j \leq n^2}.\) If the ``new" weight~\(w^{(c)}(e_j)\) is less than~\(w(e_j),\) then the corresponding minimum weight~\(M_n^{(j)}\) is also less than~\(M_n\) and the decrease is at most~\(w(e_j) - w^{(c)}(e_j).\) Moreover~\(M_n^{(j)} \neq M_n\) only if~\(e_j\) belongs to the new minimum weight planar matching~\({\cal M}^{(j)}.\) An analogous argument holds for the case~\(w(e_j) < w^{(c)}(e_j)\) and so combining we get
\[|M_n-M_n^{(j)}| \leq w(e_j) \ind\left(e_j \in {\cal M}^{(j)}\right) + w^{(c)}(e_j) \ind\left(e_j \in {\cal M}\right).\] Squaring, taking expectations and using~\((a+b)^2 \leq 2(a^2+b^2),\) we obtain
\begin{align}
\mathbb{E}\left(M_n-M_n^{(j)}\right)^2 &\leq 2\mathbb{E}w^2(e_j)\ind\left(e_j  \in {\cal M}^{(j)}\right) + 2\mathbb{E}\left(w^{(c)}(j)\right)^2\ind\left(e_j \in {\cal M}\right) \nonumber\\
&= 4\mathbb{E}w^2(e_j)\ind\left(e_j  \in {\cal M}^{(j)}\right), \label{huhu}
\end{align}
by symmetry.

By definition, we see that the minimum weight planar matching~\({\cal M}^{(j)}\) is obtained based on the edge weights~\(\{w(e_l)\}_{l \neq j} \cup \{w^{(c)}(e_j)\}\) and is therefore independent of~\(w(e_j).\) From~(\ref{huhu}) and the fact that the edge weights have uniformly bounded second moment (see Theorem statement) we therefore get that
\begin{align}
\mathbb{E}\left(M_n-M_n^{(j)}\right)^2 &\leq 4 \mathbb{E}w^2(e_j) \mathbb{P}\left(e_j \in {\cal M}^{(j)}\right) \nonumber\\
&\leq C \mathbb{P}\left(e_j \in {\cal M}^{(j)}\right) \nonumber\\
&= C \mathbb{P}\left(e_j \in {\cal M}\right), \label{huhu_two}
\end{align}
for some constant~\(C > 0\) since~\({\cal M}^{(j)}\) and~\({\cal M}\) are identically distributed.

Plugging~(\ref{huhu_two}) into~(\ref{efron_ax}), we get that
\begin{equation}\label{var_mn_ax}
var(M_n) \leq C\sum_{j=1}^{n^2} \mathbb{P}\left(e_j \in {\cal M}\right) = C \mathbb{E}N_{tot},
\end{equation}
where~\(N_{tot} \geq \tau\) is the number of edges in the minimum weight planar matching~\({\cal M}.\)  To estimate~\(N_{tot},\) we define the event
\begin{equation}\label{e_not_def}
E_{count} := \{M_n(\tau) \leq \alpha_1\tau t_n\} \bigcap \{M_n(c_n) \geq \beta_1 c_n s_n\},
\end{equation}
where~\(c_n\) is as in the statement of the Theorem. Our goal is to show that~\(N_{tot} \leq c_n+1\) with high probability and this is true if~\(c_n = n,\) since the number of edges in any planar matching of~\(K_{n,n}\) is at most~\(n.\)

Suppose next that~\(c_n = \frac{D\tau t_n}{s_n}.\) From the definition of~\(s_n\) and~\(t_n\) in~(\ref{an_bn_def_ax}), we see that~\(s_n \leq t_n\) and  we have that~\(D \geq 1\) by choice. Therefore~\(c_n \geq \tau\) and so from the deviation estimates~(\ref{dev_wn_low_ax}),~(\ref{dev_wn_up_ax}) and the union bound, we get that
\begin{equation}\label{e_not_est}
\mathbb{P}(E_{count}) \geq 1-e^{-C\tau}-e^{-2C(\log{n})^2} \geq 1-e^{-C(\log{n})^2}
\end{equation}
for some constant~\(C > 0,\) since~\(\tau \geq (\log{n})^6.\)

If~\(E_{count}\) occurs then
\begin{equation}\label{dulpo}
M_n(c_n) \geq \beta_1c_ns_n = \beta_1D\tau t_n \geq 2\alpha_1 \tau t_n \geq 2M_n(\tau) > M_n(\tau),
\end{equation}
where the first and third inequalities in~(\ref{dulpo}) are true  by the definition of~\(E_{count}\) in~(\ref{e_not_def}) and the second inequality in~(\ref{dulpo}) is true since~\(D \geq \frac{2\alpha_1}{\beta_1},\) by choice. In other words, the minimum weight of a planar matching containing at least~\(c_n\) edges is strictly larger than the minimum weight of a planar matching containing at least~\(\tau\) edges. This necessarily implies that the number of edges~\(N_{tot}\) in the minimum weight planar matching~\({\cal M}\) containing at least~\(\tau\) edges is at most~\(c_n.\) Consequently,~\[\mathbb{E}N_{tot} \leq c_n + n \mathbb{P}(E_{count}^c) \leq c_n+1\]  for all~\(n\) large, by~(\ref{e_not_est}). Plugging this into~(\ref{var_mn_ax}), we obtain the desired variance bound for~\(M_n\)  and this completes the proof of the Theorem.~\(\qed\)

Our next result concerns the saturation property described in the remark following the statement of Theorem~\ref{thm_weight_ax}. Indeed, assuming i.i.d.\ weights, we have the following convergence properties of the minimum weight of planar matchings.
\begin{theorem}\label{thm_weight_two}
Suppose the edge weights are i.i.d.\  continuous with~\(\mathbb{E}w^{2}(f) < \infty.\) For every~\(0 \leq \rho \leq 1,\) there exists a constant~\(\theta(\rho)\) such that~\(\frac{M_n(\rho n)}{n} \longrightarrow \theta(\rho)\) in   probability and in~\(L^2\) as~\(n \rightarrow \infty.\) Moreover:\\
\((a)\) The function~\(\theta(\rho)\) is non-decreasing  and convex for~\(\rho \in [0,1]\) and satisfies~\(\theta(0) = 0\) and~\(0 < \theta(\rho) < \mu \rho\) strictly, for~\(0 < \rho < 1.\)\\
\((b)\) If~\(\mathbb{E}w^{4}(f) < \infty,\) then~\(\frac{M_n(\rho n)}{n} \longrightarrow \theta(\rho)\) a.s.\ as~\(n \rightarrow \infty.\)
\end{theorem}
The property~\(\theta(\rho) < \mu \rho\) strictly, describes the ``gain" obtained by choosing minimum weight planar matchings as opposed to selecting a planar matching according to a deterministic rule.

\emph{Proof of Theorem~\ref{thm_weight_two}}\((a)\): We see that if~\(\tau = \rho n\) for some~\(0 < \rho \leq 1,\) then there are constants~\(C_1,C_2 > 0\) such that~\(C_1 \leq s_n \leq t_n \leq C_2,\) where~\(s_n\) and~\(t_n\) are as defined in~(\ref{an_bn_def_ax}). From Theorem~\ref{thm_weight_ax}, we then see that there are constants~\(D_1,D_2>0\) such that
\begin{equation}\label{sri_vaari}
D_1n \leq \mathbb{E}M_n \leq D_2n \text{ and } var(M_n) \leq D_2n.
\end{equation}
This implies that~\(\frac{M_n-\mathbb{E}M_n}{n} \longrightarrow 0\) in probability and in what follows, we use standard sub-additive arguments to obtain convergence to~\(\theta(\rho)\) and derive its properties.

Split~\(K_{n_1+n_2,n_1+n_2}\) into two disjoint subgraphs~\(K_1 \cup K_2\) with~\(K_1\) having~\(n_1\) top and bottom vertices and~\(K_2\) having~\(n_2\) top and bottom vertices. For~\(0 \leq \rho_1,\rho_2 \leq 1,\) the union of a minimum weight planar matching in~\(K_1\) containing at least~\(\rho_1n_1\) edges and a maximum weight planar matching in~\(K_2\) containing at least~\(\rho_2 n_2\) edges is a planar
matching of~\(K_{n_1+n_2,n_1+n_2}\) containing at least~\(\rho_1n_1+\rho_2n_2\) edges. Therefore
\begin{equation}\label{super_add_ax}
\mathbb{E}M_{n_1+n_2}(\rho_1n_1+\rho_2n_2) \leq \mathbb{E}M_{n_1}(\rho_1n_1) + \mathbb{E}M_{n_2}(\rho_2n_2).
\end{equation}
Setting~\(\rho_1 = \rho_2 = \rho\) and using Fekete's Lemma we get that~\(\frac{\mathbb{E}M_n}{n} \longrightarrow \theta\) for some constant~\(\theta = \theta(\rho).\) Considering the planar matching formed by the~\(\rho n\) vertical edges~\(f_i = (u_i,v_i), 1 \leq i \leq \rho n\) we get that~\(\theta(\rho) \leq \mu \rho.\) The convexity property of~\(\theta\) is obtained by setting~\(n_1 = \lambda n\) and~\(n_2 = (1-\lambda)n\) in~(\ref{super_add_ax}), dividing both sides of~(\ref{super_add_ax}) by~\(n\) and allowing~\(n \rightarrow \infty.\) Finally, since~\(M_n(\tau)\) is non-decreasing in~\(\tau,\) we get that~\(\theta\) is non-decreasing in~\(\rho.\)


\underline{\emph{Lower bound for~\(\theta(\rho)\)}}:  We first show that for every~\(\epsilon > 0\) there exists a\\\(\delta =\delta(\epsilon) > 0\) such that for any planar matching~\({\cal P}\) containing~\(m\) edges, we have
\begin{equation}\label{t_pi_ax}
\mathbb{P}\left(W({\cal P}) \leq \delta m\right) \leq e^{-\epsilon m},
\end{equation}
where we recall that~\(W({\cal P})\) is the weight of~\({\cal P}.\)

Indeed, let~\({\cal P} = \{e_1,\ldots,e_m\}\) so that~\(W({\cal P}) = \sum_{i=1}^{m} w(e_i).\) Using the Chernoff bound we obtain for~\(\delta,s > 0\) that
\begin{equation}
\mathbb{P}(W({\cal P}) \leq \delta m)  = \mathbb{P}\left(\sum_{i=1}^{m}w(e_i) \leq \delta m\right) \leq e^{s\delta m}\prod_{i=1}^{m}\mathbb{E}\left(e^{-sw(e_i)}\right).\label{y_1_eq1}
\end{equation}
For a fixed \(\eta > 0\) to be determined later, we write~\[\mathbb{E}e^{-sw(e_i)} = \int_{w(e_i) < \eta} e^{-sw(e_i)} d\mathbb{P} + \int_{w(e_i) \geq \eta} e^{-sw(e_i)} d\mathbb{P} \]
and use~\[\int_{w(e_i) < \eta} e^{-sw(e_i)} d\mathbb{P} \leq \mathbb{P}(w(e_i) < \eta) \text{ and }\int_{w(e_i) \geq \eta} e^{-sw(e_i)} d\mathbb{P} \leq e^{-s\eta} \]
to get that
\begin{equation}\label{hannah}
\mathbb{E}e^{-sw(e_i)} \leq \mathbb{P}(w(e_i) < \eta) + e^{-s\eta}.
\end{equation}

Since the weights are continuous, we can choose~\(\eta > 0\) small so that~\[\mathbb{P}(w(e_i) < \eta) \leq \frac{e^{-6\epsilon}}{2}.\] Fixing such an~\(\eta\) we choose~\(s = s(\eta,\epsilon) > 0\) large so that the second term in~(\ref{hannah})~\(e^{-s\eta} < \frac{e^{-6\epsilon}}{2}.\) Combining we then get from~(\ref{hannah}) that~\(\mathbb{E}e^{-st(e_i)} \leq e^{-6\epsilon}\) and so from~(\ref{y_1_eq1}) we obtain~\[\mathbb{P}(W({\cal P}) \leq \delta m) \leq e^{s\delta m} e^{-6\epsilon m} \leq e^{-2\epsilon m}\] for all \(m \geq 1,\) provided \(\delta = \delta(s,\epsilon) > 0\) is small. This proves~(\ref{t_pi_ax}).

For any planar matching~\({\cal P}\) containing at least~\(\tau \geq \rho n\) edges, we now set~\(\epsilon = \frac{4}{\rho}\) in~(\ref{t_pi_ax}) to get that~\[\mathbb{P}(W({\cal P}) \leq \delta n) \leq \exp\left(-\frac{4\tau}{\rho}\right) \leq e^{-4n}\] for some~\(\delta = \delta(\rho) > 0.\) There are at most~\({n \choose \tau} \leq 2^{n}\) choices for the top vertices of~\({\cal P}\) and at most~\({n \choose \tau} \leq 2^{n}\) choices for the bottom vertices of~\({\cal P}.\) Therefore if~\(E_{bad}\) is the event that there exists a planar matching containing~\(\tau \geq \rho n\) edges and having weight at most~\(\delta n,\) then
\begin{equation}\label{e_bad_est}
\mathbb{P}(E_{bad}) \leq \sum_{\rho n \leq \tau \leq n} \left(2^{n}\right)^2 e^{-4n} \leq n4^{n}e^{-4n} \leq e^{-n}.
\end{equation}

From~(\ref{e_bad_est}), we see that~\[\mathbb{E}M_n(\rho n) \geq \mathbb{E}M_n(\rho n) \ind(E^c_{bad}) \geq \delta n (1-e^{-n}) \geq \frac{\delta}{2} n\] and so~\(\theta(\rho) \geq \frac{\delta}{2} > 0.\) This completes the proof of the lower bound for~\(\theta(\rho).\)

\underline{\emph{Upper bound for~\(\theta(\rho)\)}}: We are to show that~\(\theta(\rho) < \mu \rho\) strictly for each\\\(0 < \rho < 1\) and we proceed as follows. Let~\(M_n = M_n(\rho n)\) be the minimum weight of a planar matching containing at least~\(\rho n\) edges and let~\(r\) be the integer satisfying~\(\frac{1}{r} \leq \rho < \frac{1}{r-1}.\)  From the super-additivity relation~(\ref{super_add_ax}), we get that~\(\mathbb{E}M_{rn} \leq n \mathbb{E}M_r\)
and so
\[\theta(\rho) = \lim_{n} \frac{\mathbb{E}M_{rn}}{rn}  \leq \frac{\mathbb{E}M_r}{r}.\] Therefore to show that~\(\theta(\rho) < \mu \rho\) strictly, it is enough to show that~\(\mathbb{E}M_r < \mu\) strictly, since this would imply that~\(\theta(\rho) < \frac{\mu}{r} \leq \mu \rho.\)

By definition~\(1 \leq \rho r < 2\) strictly and so~\(M_r\) is the minimum weight of an edge in the bipartite graph~\(K_{r,r}.\) There are~\(r^2\) edges in~\(K_{r,r}\) and each edge~\(f \in K_{r,r}\) is the minimum weight edge with probability~\(\frac{1}{r^2}.\) Therefore
\begin{equation}\label{dodo}
\mathbb{E}M_r = \sum_{f }\mathbb{E}w(f) \ind(E_{low}(f))
\end{equation}
where~\(E_{low}(f)\) is the event that~\(f\) is the minimum weight edge in~\(K_{r,r}\) and the summation is over all edges in~\(K_{r,r}.\) If~\(E_{low}(f)\) occurs, then~\(w(f) < \frac{1}{r^2} \sum_{h} w(h)\) strictly since the weights are continuous and so from~(\ref{dodo}) we get that
\begin{align}
\mathbb{E}M_r &< \sum_{f}\mathbb{E}\frac{1}{r^2}\sum_{h} w(h) \ind(E_{low}(f)) \nonumber\\
&= \frac{1}{r^2} \mathbb{E}\sum_{h} w(h) \sum_{f} \ind(E_{low}(f)) \nonumber\\
&= \frac{1}{r^2} \mathbb{E} \sum_{h} w(h) \nonumber\\
&= \mu.
\end{align}
This obtains the desired upper bound for~\(\theta(\rho)\) and completes the proof of part~\((a)\) of the Theorem.~\(\qed\)

\emph{Proof of Theorem~\ref{thm_weight_two}\((b)\)}: Recalling that~\(\{e_i\}_{1 \leq i \leq n^2}\) is the set of~\(n^2\) edges in the complete bipartite graph~\(K_{n,n},\)  we say that~\(e_j\) is \emph{pivotal} if every minimum weight planar matching of~\(K_{n,n}\) contains~\(e_j.\) Letting~\(L_n\) denote the number of pivotal edges in~\(K_{n,n},\) we have from Lemma~\(5.4\) of~\cite{chat} that
\begin{equation}\label{emily}
\mathbb{E}(M_n-\mathbb{E}M_n)^4 \leq C_1 \mathbb{E}L^2_n\left(\mathbb{E}w^2(e_1)\right)^2 + C_1 \mathbb{E}L_n\mathbb{E}w^4(e_1).
\end{equation}
Using~\(L_n \leq n,\) the maximum size of a planar matching and the fact that edge weights have bounded fourth moments, we get from~(\ref{emily}) that
\begin{equation}\label{mn_mom_four}
\mathbb{E}(M_n-\mathbb{E}M_n)^4 \leq C_2n^2
\end{equation}
for some constant~\(C_2 > 0.\) From the Borel-Cantelli lemma, we then get that~\[\frac{M_n-\mathbb{E}M_n}{n} \longrightarrow 0 \text{ a.s. }\] as~\(n \rightarrow \infty\) and so the desired almost sure convergence in the statement of the Theorem then follows from the definition of~\(\theta(\rho).\) This completes the proof of the Theorem.~\(\qed\)

\setcounter{equation}{0}
\renewcommand\theequation{A.\arabic{equation}}
\section*{Appendix}
\emph{Proof of~(\ref{mn_est})}: To obtain the deviation upper bound in~(\ref{mn_est}), we use a direct counting argument. Let~\(p = p_{up}\) and suppose that~\(\{M_n \geq t\}\) occurs for some~\(t \geq 4en\sqrt{p}.\) Let~\({\cal S} = \{(u_i,v_i)\}_{1 \leq i \leq t}, t \geq 4en\sqrt{p}\) be an~\(L-\)constrained planar matching. The number of choices for~\(\{v_i\}_{1 \leq i \leq t}\) is~\({n \choose t}\) and given~\(\{v_i\}_{1 \leq i \leq t},\) the number of choices of~\(\{u_i\}_{1 \leq i \leq t}\) is at most~\({n \choose t}.\)  Therefore the total number of choices for~\({\cal S}\) is at most~\({n \choose t}^2 \) and since each edge of~\(K_{n,n}\) is present with probability~\(p\) independently of other edges, we use the estimate~\({n \choose k} \leq \left(\frac{ne}{k}\right)^{k}\) to get that
\begin{align}\label{up_bound2}
\mathbb{P}(R_n \geq t) &\leq  {n \choose t} \cdot {n \choose t} \cdot p^{t} \nonumber\\
&\leq \left(\frac{ne}{t}\right)^{2t} p^{t} \nonumber\\
&= \left(\frac{ne\sqrt{p}}{t}\right)^{2t} \nonumber\\
&\leq \left(\frac{1}{2}\right)^{2t}.
\end{align}
This proves the deviation bound in~(\ref{mn_est}).

To obtain the upper bound for~\(\mathbb{E}R_n,\) we proceed as follows. If~\(4en\sqrt{p} \geq \frac{1}{2},\) then we use~(\ref{up_bound2}) and the inequality~\(\mathbb{E}M_n \leq 1+\sum_{t \geq 1} \mathbb{P}(M_n \geq t)\) to get that
\begin{align}
\mathbb{E}R_n &\leq 1+4en\sqrt{p} + \sum_{t \geq 4en\sqrt{p}+1} \frac{1}{2^{t}} \nonumber\\
&\leq 1+4en\sqrt{p} + \frac{1}{2^{4en\sqrt{p}}} \nonumber\\
&\leq 4en\sqrt{p} + 2. \nonumber
\end{align}
Since~\(4en\sqrt{p} \geq \frac{1}{2},\) we get that~\(\mathbb{E}R_n \leq 16en\sqrt{p}.\) If on the other hand~\(4en\sqrt{p} \leq \frac{1}{2},\) then from~(\ref{up_bound2}) we get that~\(\mathbb{P}(R_n \geq t) \leq \left(4en\sqrt{p}\right)^{t} \) for any~\(t \geq 1\) and so using~\[\mathbb{E}R_n = \sum_{t \geq 1} t \mathbb{P}(R_n = t) \leq \sum_{t \geq 1} t \mathbb{P}(R_n \geq t)\]
we get that~\[\mathbb{E}R_n \leq \sum_{t \geq 1} t(4en\sqrt{p})^{t}  = \frac{4en\sqrt{p}}{(1-4en\sqrt{p})^2} \leq 16en\sqrt{p}.\] This obtains the expectation upper bound in~(\ref{mn_est}) and therefore completes the proof of~(\ref{mn_est}).~\(\qed\)

\subsection*{\em Acknowledgement}
I thank Professors Rahul Roy, Federico Camia, Alberto Gandolfi, C. R. Subramanian and the referees for crucial comments that led to an improvement of the paper. I also thank IMSc and IISER Bhopal for my fellowships.

\subsection*{\em Data Availability Statement}
Data sharing not applicable to this article as no datasets were generated or analysed during the current study.

\subsection*{\em Conflict of Interest}
The authors have no conflicts of interest to declare that are relevant to the content of this article. No funding was received to assist with the preparation of this manuscript.

\bibliographystyle{plain}

\end{document}